\documentclass[11pt]{article}
\usepackage{anysize}
\marginsize{2.0cm}{2.0cm}{2.0 cm}{2.0 cm}
\usepackage{setspace}

\usepackage{latexsym,amsfonts,amsmath}
\usepackage{graphicx}
\usepackage{color}
\usepackage{epstopdf}
\usepackage[titletoc]{appendix}

\def\RR{\mathbb{R}}

\def\liminf{\mathop{\underline{\lim}}}

{\bfseries}{\itshape}

{\bfseries}{\itshape}

\newtheorem{theorem}{Theorem}[section]{\bfseries}{\itshape}

{\bfseries}{\itshape}

\newtheorem{proposition}{Proposition}[section]{\bfseries}{\itshape}

\newtheorem{example}{Example}[section]{\bfseries}{\itshape}

\newtheorem{lemma}{Lemma}[section]{\bfseries}{\itshape}

{\bfseries}{\itshape}

\newtheorem{definition}{Definition}[section]{\bfseries}{\itshape}

{\bfseries}{\itshape}

\begin{document}
\title{On the structure of optimal solutions in a mathematical programming problem in a convex space}

\author{Alexey Piunovskiy\thanks{Department of Mathematical Sciences, University of
Liverpool, Liverpool, U.K.. E-mail: piunov@liv.ac.uk.}~ and Yi
Zhang \thanks{School of Mathematics, University of
Birmingham, Birmingham, B15 2TT, U.K.. E-mail: y.zhang.29@bham.ac.uk.}}
\date{}
\maketitle

\par\noindent{\bf Abstract:} We consider an optimization problem in a convex space $E$ with an affine objective function, subject to $J$ constraints in the forms of inequalities on some other affine functions, where $J$ is a given nonnegative integer. Under suitable conditions, we apply the Feinberg-Shwartz lemma in finite dimensional convex analysis to show that there exists an optimal solution, which is in the form of a mixture of no more than $J+1$ extreme points of $E$. It seems that in the current setup, this result has not yet been made available, because the concerned problem does not fit into the framework of standard convex optimization problems.
\bigskip

\bigskip


\section{Introduction}\label{PZ23aSecIntroduction}

In the literature of convex analysis and optimization, see e.g., \cite{Aliprantis:2007, andenash, luen, Rockafellar:1970, Rockafellar:1974}, the standard definition of a convex set is a convex subset of a vector space. It has been noticed that the requirement that a convex set to be embedded into a vector space is restrictive for several applications in e.g., decision theory and petroleum engineering, see \cite{Gudder:1979,Gudder:1980} for some discussions and historical remarks. The natural generalization is to consider the so called convex space or convex structure, where only operations of finite convex combinations are defined, without involving a vector space, see \cite{Fritz:2015}. Affine functions can be defined on a convex space in the same way as on a convex set in a vector space.  According to \cite[Theorem 4]{Gudder:1977}, a convex space is isomorphic to a convex set in a cone, rather than a vector space, where the cone is defined again without referring to a vector space, see precise definitions in Subsection \ref{PZ23aSubSec01}. Therefore, in terms of terminology, we will talk about convex sets in a cone, instead of convex spaces.

In this paper, we consider an optimization problem in a convex set of a cone with an $(-\infty,\infty]$-valued affine objective function on it, and subject to constraint inequalities on other $J$ affine functions, where $J\ge 0$ is a fixed nonnegative integer.  

Our contributions are as follows. For an optimization problem in a convex set $E$ of a cone with an $(-\infty,\infty]$-valued affine objective function on it, and subject to constraint inequalities on other $J$ affine functions, we impose suitable conditions under which we show that this problem has an optimal solution which is in the form of the mixture of no more than $J+1$ extreme points in $E$.  In \cite[Subsection 3.2.5]{APZY}, the similar result was obtained for $E$ being a metrizable compact convex subspace of a topological Hausdorff space, where the latter requirement was inaccurately missing.  Here, the conditions are more general in that we assume the consistency of the problem, and the compactness of the convex set $E$ in a cone, which is not necessarily metrizable, as well as the totality of the space of bounded below lower semicontinuous affine functions on $E$. The application of this result to Markov decision processes will be demonstrated in a companion paper. This paper can also be viewed as an application of the Feinberg-Shwartz lemma  in \cite{FeinbergS:1996} in finite dimensional convex analysis, quoted as Proposition \ref{lb31} below, on which our proof is based.

The rest of the paper is organized as follows. In Section \ref{PZ23aSec01}, we provide basic definitions and present the main result. Section \ref{PZ23aSec02} collects auxiliary lemmas and facts. The proof of the main theorem is presented in Section \ref{PZ23aSec05}. The paper ends with a conclusion in Section \ref{PZ23aSec06}.

\section{Problem statement and main result}\label{PZ23aSec01}

In this section, we state the concerned optimization problem. Before that, we provide some definitions that are directly involved in the problem statement and main result to be presented.  Although the objects to be defined  are all standard in the context of convex sets in a vector space, we need to consider their natural generalizations in the context of convex sets in a cone, or say convex spaces.

\subsection{Basic definitions}\label{PZ23aSubSec01}

To begin with, we define a cone without involving a vector space as follows.

\begin{definition}[Cone]\label{db311}
Consider a nonempty set $X$ with an element ${\bf 0}\in X$ and the operations of addition, denoted as $x+y\in X$ for $x,y\in X$, and scalar multiplication by non-negative real numbers, denoted as $ax\in X$ for all $x\in X$ and $a\in  [0,\infty)$. The set $X$ is called a cone (over $[0,\infty)$) if it satisfies the following conditions:
\begin{itemize}
\item $a(x+y)=ax+ay$ for all $a\in[0,\infty)$ and $x,y\in X$;
\item $a(bx)=(ab)x$ for all $a,b\in[0,\infty)$ and  $x\in X$;
\item $1x=x$ and  $0x={\bf 0}$ for all $x\in X$;
\item $x+y=y+x$ for all $x,y\in X$;
\item $x+(y+z)=(x+y)+z$ for all $x,y,z \in X$;
\item $(a+b)x=ax+bx$ for all $a,b\in [0,\infty)$ and $x\in X$;
\item $x+{\bf 0}=x$ for all $x\in X.$
\end{itemize}
\end{definition}

The above definition is close to the one of a cone with a zero on p.93 of \cite{Gudder:1979}. The only difference is that in  \cite{Gudder:1979}, a cone is further required to be cancellative.  A cone is called cancellative if for any $y,z\in X$,  if for some $x\in X,$ it holds that $x+y=x+z$, then $y=z.$ Here we do not require this property to hold for the concerned set to be a cone.

An example of a cone (over $[0,\infty)$) which is not cancellative is $X=[0,\infty]$ with the usual addition and multiplication, ${\bf 0}=0$,  $0\cdot \infty:=0$,  because for $y,z\in X,$ $\infty+y=\infty+z$ (both sides being $\infty$) does not imply $y=z.$

\begin{definition}[Convex set in a cone]\label{db312}
A subset $E$ of a cone $X$ is called a convex set (in $X$) if $\alpha x+(1-\alpha) y\in E$ for all $x,y\in E$ and $\alpha\in[0,1]$. A subset of $E$ is called a convex subset of $E$ if it is a convex set in $X$.
\end{definition}

Let us say a few words about terminologies. In view of \cite[Theorem 4]{Gudder:1979}, we may well call the above defined convex set in a cone a convex space: this is consistent with the definition of a convex space in \cite{Fritz:2015}.
In the literature, by a convex set it is often meant a convex subset of a vector space. A vector space is a particular cone. Here and below, by a convex set, unless stated otherwise, we always mean a convex set in a cone, not necessarily in a vector space. Some results need the convex set to be a subset of a vector space, see, for instance,  Example \ref{TeenfestivalExample01} below. When this is to be emphasized, we will say explicitly that the convex set is in the underlying vector space.

For example,  the set $X$ of all $[0,\infty]$-valued measures on a given measurable space $(\Omega,{\cal F})$, with the usual addition and scalar multiplication, is a non-cancellative cone, and one can consider convex subsets of this convex cone of $[0,\infty]$-valued measures on a $\sigma$-algebra $\cal B$. This is the situation that arises in studies of Markov decision processes with total undiscounted criteria, where the occupation measures form convex set in the cone of $[0,\infty]$-valued measures, which in general cannot be embedded into a vector space, see \cite{Dufour:2020,Dufour:2012}.

\begin{definition}[Face and extreme point]\label{db32}
Let $E$ be a nonempty convex set. A nonempty subset $D\subseteq E$ is called an extreme subset of $E$, if for each point $u\in D$, the representation $u=\alpha u_1+(1-\alpha) u_2$ with $\alpha\in(0,1)$, $u_1,u_2\in E$ implies that $u_1,u_2\in D.$ A point $u\in E$ is called an extreme point of $E$ if the singleton $\{u\}$ is an extreme subset  of $E$. An extreme subset of $E$ may be not convex. E.g., a pair of distinct extreme points forms such a non-convex extreme subset. A convex extreme subset of $E$ is called a face of $E$.
\end{definition}
Suppose in the above definition, the convex set $E$ has a topology on it. Then by a closed extreme subset (or face) of $E$ we mean a closed subset of $E$, which is also an extreme subset (respectively, a face) of $E$.

\begin{definition}[Affine function]\label{db34}
A $(-\infty,+\infty]$-valued function $f(\cdot)$ defined on a convex set $E$ in a cone is called affine if
$$f(\alpha x_1+(1-\alpha) x_2)=\alpha f(x_1)+(1-\alpha)f(x_2)$$
for all $x_1,x_2\in E$ and all $\alpha\in [0,1]$.
\end{definition}
In the above definition, we used the same notation for the addition and scalar multiplication in the underlying cone and in $(-\infty,\infty]$.  The context should exclude any confusion.

\subsection{Problem statement and main result}
Now we state the concerned optimization problem.

Consider a nonempty convex set $E$ in a cone, and  endow $E$ with some topology. This topology on $E$ needs not be consistent with the linear operations (addition and scalar multiplication).
Let  $\hat{\cal C}(E)$ be the family of $(-\infty, +\infty]$-valued bounded from below lower semicontinuous affine functions on $E$. As in \cite[p.222]{Gudder:1977}, we call  $\hat{\cal C}(E)$ total if for each $x,y\in E$ with $x\ne y$, there is some function $f\in \hat{\cal C}(E)$  such that $f(x_1)\ne f(x_2).$

Let $W_0(\cdot),W_1(\cdot),\ldots,W_J(\cdot)\in\hat{\cal C}(E)$ be given. We consider the following optimization problem
\begin{eqnarray}\label{eB12p}
&&\mbox{Minimize over $x\in E$: } W_0(x) \nonumber\\
&&\mbox{ subject to } W_j(x)\le d_j,~j=1,2,\ldots, J,
\end{eqnarray}
where $d_j\in\RR$ are fixed constants and $J\ge 0$. The case of $J=0$ corresponds to the absence of constraint inequalities in problem (\ref{eB12p}).

In the context of problem (\ref{eB12p}), any point in $E$ will be referred to as a solution (for  problem (\ref{eB12p})).  A solution is called feasible for problem (\ref{eB12p}) if it satisfies the constraints in it. If $J=0$, then every solution is feasible. Problem (\ref{eB12p}) is called consistent, if it has some feasible solutions. A feasible solution $x^\ast\in E$ is called an optimal solution for problem (\ref{eB12p}) if $W_0(x^\ast)\le W_0(x)$ for all feasible solutions $x\in E$.

Now we present our main result.

\begin{theorem}\label{tB1}
Suppose that the nonempty convex set $E$  (in a cone) is compact,  problem (\ref{eB12p}) is consistent, and $\hat{\cal C}(E)$ is total. Then there exists an optimal solution, say $x^\ast\in E$, for  problem (\ref{eB12p}) in the form of a convex combination of at most $J+1$ extreme points in $E$, i.e., $x^\ast=\sum_{k=1}^{J+1} \alpha_k x_k$, where $\alpha_k\in[0,1]$, $\sum_{k=1}^{J+1} \alpha_k=1$, and $x_k$ is extreme in $E$ for each $k=1,2,\ldots,J+1$.
\end{theorem}

\section{Preliminary results}\label{PZ23aSec02}

In this section, we present some lemmas and recall some further facts, which are needed in the proof of the main theorem. Lemmas \ref{PZ23aRem01} and \ref{PZ23aLemma01} are immediate consequences of the relevant definitions. However, we still provide their complete proofs in hope of further improving the readability.

\begin{lemma}\label{PZ23aRem01}
If $E_1$ is an extreme subset of $E_0$, whereas $E_0$ is an extreme subset of a nonempty convex set $E$, then $E_1$ is also an extreme subset of $E$. In particular, a face of a face of a nonempty convex set $E$ is also a face of $E$.
\end{lemma}
\par\noindent\textit{Proof}. We verify the first assertion, only, as the second assertion follows from it.
Let $x\in E_1$ and assume that $x=\alpha x_1+(1-\alpha)x_2$ with $\alpha\in(0,1)$ and $x_1,x_2\in E$. Then  $x_1,x_2\in E_0$ because $x\in E_0$, and $E_0$ is an extreme subset of $E$. This in turn leads to $x_1,x_2\in E_1$ because $E_1$ is an extreme subset of $E_0$. Thus $E_1$ is an extreme subset of $E$, as required.  $\hfill\Box$

\begin{lemma}\label{PZ23aLemma01}
Consider a nonempty convex set $E$, which is endowed with a topology. Let $I$ be an arbitrary nonempty index set. If $E_i$, $i\in I$, are extreme subsets (or faces) of $E$ and $\bigcap_{i\in I} E_i\ne \emptyset$, then $\bigcap_{i\in I} E_i$ is an extreme subset (respectively, face) of $E$.  In particular, if $E_i$, $i\in I$, are closed extreme subsets (or faces) of $E$ and $\bigcap_{i\in I} E_i\ne \emptyset$, then $\bigcap_{i\in I} E_i$ is a closed extreme subset (respectively, face) of $E$.
\end{lemma}
\par\noindent\textit{Proof.} Suppose $x\in \bigcap_{i\in I} E_i$, $x=\alpha x_1+(1-\alpha) x_2$ with $\alpha\in(0,1)$, and $x_1,x_2\in E$.  Then, for each $i\in I$, from the inclusion $x\in E_i$, we see that  $x_1,x_2\in E_i$, because $E_i$ is an extreme subset of $E$. Since the inclusion $x\in E_i$ holds for all $i\in I,$ we see that $x_1,x_2\in \bigcap_{i\in I} E_i$. Thus, $\bigcap_{i\in I} E_i$ is an extreme subset of $E$.

If $E_i$ are all convex, then so is their intersection $\bigcap_{i\in I} E_i$, and in that case,  $\bigcap_{i\in I} E_i$ will be a face of $E$.

The last asserton holds because the set $\bigcap_{i\in I} E_i$ is closed if the sets $E_i$ are all closed.
$\hfill\Box$

Given a nonempty convex set $E$ (in a cone $X$), the minimal face of $E$ that contains a point $u\in E$, i.e., the intersection of all the faces of $E$ containing $u$, is denoted by $G_E(u)$.

\begin{definition}[Pareto optimality]\label{db31}
Let $E$ be a fixed nonempty convex set in the cone $(-\infty,\infty]^{J+1}$ with some nonnegative integer $J$. A point $u\in E$ is called Pareto optimal if, for each $v\in E$, the componentwise inequality $v\le u$ implies that $v=u$. The collection of all Pareto optimal points is denoted by $Par(E)$.
\end{definition}



Let $J\ge 0$ be a fixed nonnegative integer.  Consider a nonempty convex  $E$ in $\RR^{J+1}$, and $\vec{u}=(u_0,\dots,u_J)\in Par(E)$. The next useful result of Feinberg and Shwartz, see Lemma 3.2 of \cite{FeinbergS:1996}, gives a structure of $G_E(\vec{u}).$

\begin{proposition}[Feinberg and Shwartz]\label{lb31}
Let $J\ge 0$ be a nonnegative integer. Suppose $E$ is a fixed nonempty convex subset of $\RR^{J+1}$, and $\vec{u}=(u_0,\dots,u_J)\in Par(E)$. Then the following assertions are valid.
\begin{itemize}
\item[(a)] $G_E(\vec{u})\subseteq Par(E)$.
\item[(b)] For some $1\le k\le J+1$, and $1\le i\le k$, there are some $\vec{b}^i=(b^i_0,\dots,b^i_J)\in \RR^{J+1}$ and $\beta^i\in \mathbb{R}$, defining the hyperplanes in $\RR^{J+1}$ \begin{eqnarray*}
H^i=\left\{\vec{x}=(x_0,\dots,x_J)\in\mathbb{R}^{j+1}:~\sum_{j=0}^J b^i_j x_j =\beta^i\right\},~i=1,2,\ldots, k,
\end{eqnarray*}
satisfying the following properties:
\begin{itemize}
\item[(i)] $\vec{b}^i\ge 0$, $\sum_{j=0}^J b^i_j>0$, for $i=1,2,\ldots, k-1$ and $\vec{b}^k>0$. Here all the inequalities are componentwise.
\item[(ii)] $\sum_{j=0}^J b_j^1 x_j\ge \beta^1$ for all $\vec{x}=(x_0,\dots,x_J)\in E^0:=E$, $\sum_{j=0}^J b_j^1 u_j=\beta^1$; for $i=1,2,\ldots,k-1$,   $\sum_{j=0}^J  b_j^{i+1} x_j\ge \beta^{i+1}$ for all $\vec{x}=(x_0,\dots,x_J)\in E^i:=E^{i-1}\cap H^i$, and $\sum_{j=0}^J b_j^{i+1} u_j=\beta^{i+1}$.
\item[(iii)] $G_E(\vec{u})=E^k:=E^{k-1}\cap H^k$.
\end{itemize}
\end{itemize}
\end{proposition}
\par\noindent\textit{Proof.}  See Lemmas 3.1 and 3.2 of \cite{FeinbergS:1996}. $\hfill\Box$

The proof of Theorem \ref{tB1} is based on an application of the result of Feinberg and Shwartz quoted in Proposition \ref{lb31}. We shall refer to it as the Feinberg-Shwartz lemma.

The correctness of the previous result requires $E$ to be a convex subset of $\mathbb{R}^n.$ More precisely, the next example shows that if $E$ is a nonempty convex subset of the cone $(-\infty,\infty]^2$ and $\vec{u}\in \mathbb{R}^2$, then Proposition \ref{lb31}  may fail to hold.
\begin{example}\label{TeenfestivalExample01}
Let $E'$ be a closed disk of radius $0.5$ and centered at $(1,1.5)$, and ${E}=E'\cup \{(\infty,v_2): v_2\in[0,2]\}$. Then this set $E$ is convex in the cone $(-\infty,\infty]^2$, and $\vec{u}=(1,1)$ is Pareto in $E$. It is also an extreme point of $E$. Thus, $G_E(\vec{u})=\{\vec{u}\}.$

If Proposition \ref{lb31} was applicable to this convex set $E$ in the cone $(-\infty,\infty]^2$ and the point $\vec{u}=(1,1)$, then there would be nonnegative $b_1,b_2\ge 0$ with $b_1+b_2>0$ and a constant $\beta$ such that
$b_1 \cdot 1 +b_2\cdot 1=\beta$ and
$b_1 v_1+b_2 v_2\ge \beta$ for all $\vec{v}=(v_1,v_2)\in E$.

The requirement of $b_1 \cdot 1 +b_2\cdot 1=\beta$ and $b_1 v_1+b_2 v_2\ge \beta$ for all $\vec{v}=(v_1,v_2)\in E'\subset E$ imply that $b_1=0,$ and $b_2=\beta>0$, as given by the unique supporting hyperplane in $\mathbb{R}^2$ of the disk $E'$ at $\vec{u}=(1,1).$


On the other hand, since $0\cdot\infty:=0,$ for $(\infty,0)\in E$,  $b_1\cdot \infty+b_2 \cdot 0=0<\beta,$ yielding a contradiction.
\end{example}

The next result, in particular, asserts that a nonempty convex compact set $E$ in a cone always has some extreme points, provided that $\hat{\cal C}(E)$ is total. Its proof is similar to the one of \cite[Lemma 7.65]{Aliprantis:2007}.
\begin{lemma}\label{lB6}
Suppose $E$ is nonempty, compact and convex set in a cone. Recall that  $\hat{\cal C}(E)$ is the space of $(-\infty,+\infty]$-valued bounded from below lower semicontinuous affine functions on $E$. Then the following assertions are valid.
\begin{itemize}
\item[(a)] Suppose $F$ is a closed extreme subset (or a closed face) of $E$ and $f(\cdot)\in\hat{\cal C}(E)$. Then
\begin{eqnarray*}
Y:=\left\{x\in F:~ f(x)=\inf_{\tilde x\in F} f(\tilde x)\right\}
\end{eqnarray*}
is a closed extreme subset (respectively, a closed face) of $E$. 
\item[(b)] Assume that $\hat{\cal C}(E)$ is total, i.e., if $x_1\ne x_2\in E$, then there is $f(\cdot)\in\hat{\cal C}(E)$ such that $f(x_1)\ne f(x_2)$. Then each closed face  $F$ of $E$ contains at least one extreme point $\hat x$ of $E$, and the singleton $\{\hat x\}$ is closed. In particular, $E$ has some extreme points.
\end{itemize}
\end{lemma}

\par\noindent\textit{Proof.} (a) Let $\inf_{\tilde x\in F} f(\tilde x):=a\in (-\infty,\infty].$  If $a=+\infty$, then $Y=F$ and the statement follows from the assumed properties of $F$.

Below, we assume that $a\in\RR$. Then $Y:=\left\{x\in F:~ f(x)=a\right\}=\left\{x\in F:~ f(x)\le a\right\}$.

Firstly, consider $F$ as a closed extreme subset of $E$.

It follows that the set $Y$ is nonempty and closed, because $F$ is closed, and the function $f(\cdot)$ is lower semicontinuous.

Let us verify that $Y$ is an extreme subset of $E$. Indeed, if $x\in Y$ and $x=\alpha x_1+(1-\alpha)x_2$ with $\alpha\in(0,1)$ and $x_1,x_2\in F$, then, $a=f(x)=\alpha f(x_1)+(1-\alpha)f(x_2)$ because $f(\cdot)$ is affine on $E$.  Since for each $\tilde x\in F$, $f(\tilde x)\ge a$, to have $x\in Y$ (i.e., $f(x)=a$), we must have $f(x_1)=f(x_2)=a$	 leading to $x_1,x_2\in Y$.  This shows that $Y$ is an extreme subset of $F$. Now, according to Lemma \ref{PZ23aRem01}, $Y$ is also an extreme subset of $E$.

Thus, $Y$ is a closed extreme subset of $E.$

Finally, we note that  if $F$ is a closed face of $E$, then the set $Y$ is convex. Indeed, if $x_1,x_2\in Y$, then  for each $\alpha\in [0,1]$, $\alpha x_1+(1-\alpha)x_2\in F$, and $f(\alpha x_1+(1-\alpha x_2))=\alpha f(x_1)+(1-\alpha)f(x_2)=a$, where the first equality holds because $f(\cdot)$ is affine, and the last equality holds because $x_1,x_2\in Y$. Consequently, $\alpha x_1+(1-\alpha )x_2\in Y$.  It follows from this and what was established earlier that $Y$ is a closed face of $E.$

(b) Fix a closed face $F$ of $E.$ Consider the family ${\cal A}$ of closed faces of $F$, which is partially ordered with respect to the inclusion, i.e., $A\preceq B$ if and only if $B\subseteq A$. For each chain in ${\cal A}$, since $F$ is compact (as a closed subset of the compact $E$) and each finite subfamily say $A_1\supseteq\dots\supseteq A_n$ of the chain in ${\cal A}$ has a nonempty intersection $A_n\ne\emptyset$, one may refer to
Theorem 2.31 of \cite{Aliprantis:2007} for that each chain in ${\cal A}$ has a nonempty intersection. According to Lemma \ref{PZ23aLemma01}, this intersection is an upper bound of the chain in ${\cal A}$.  Therefore, according to Zorn's Lemma, see e.g., \cite[Lemma 1.7]{Aliprantis:2007}, ${\cal A}$ has a maximal element, say $\hat{F}$, and no proper subset of $\hat{F}$ can be a closed face of $F$.

Next, we show that $\hat{F}$ is a singleton. Suppose for contradiction that there exist $x_1, x_2\in\hat F$ such that $x_1\ne x_2.$ Since by assumption $\hat{\cal C}(E)$ is total, we may consider some $f(\cdot)\in\hat{\cal C}(E)$ such that $f(x_1)\ne f(x_2)$. Without loss of generality, assume that $f(x_1)>f(x_2)$. Now by part (a), applied to the nonempty compact and convex set $F$, we see that
\begin{eqnarray*}
Z:=\left\{ x\in\hat F:~f(x)=\inf_{\tilde x\in\hat F} f(\tilde x)\le f(x_2)\right\}
\end{eqnarray*}
is a closed face of $F$, and $x_1\notin Z$. As a result,
$Z\subset \hat F$ is a proper subset of $\hat F$ as well as a closed extreme subset of $F$, which is a desired contradiction. We conclude that $\hat F=\{\hat x\}$ is a singleton, closed in $F$ and thus in $E$, and the point $\hat x\in F$ is an extreme point of $F$. By Lemma \ref{PZ23aRem01}, $\hat{x}$ is also an extreme point of $E$.
 $\hfill\Box$

It is known, see e.g., \cite[Corollary 7.66]{Aliprantis:2007} that every nonempty convex compact subspace $E$ of a locally convex Hausdorff space has some extreme points.  Lemma \ref{lB6} is an extension of that result, because in a locally convex Hausdorff space, the set of real-valued continuous affine functions separate points: see Item 5 of \cite[Ch.XI,\S1]{meyer} or \cite[Corollary 5.82]{Aliprantis:2007}.

\section{Proof of Theorem \ref{tB1}}\label{PZ23aSec05}
In this section, we provide the detailed proof of Theorem  \ref{tB1}.

\par\noindent\textit{Proof of Theorem \ref{tB1}.}  First, we assume that there is a feasible solution, say $\underline{x}\in E$, with  $W_0(\underline{x})<\infty$. Then there is an optimal solution, say $\overline{x}\in E$, for problem  (\ref{eB12p})  with a finite value, say $W_0(\overline{x})=:d_0\in \mathbb{R}.$ It will be explained at the end of this proof that this assumption can be withdrawn.

Introduce the space of performance vectors
\begin{eqnarray*}
{\bf O}:=\{\vec W(x)=(W_0(x),W_1(x),\ldots,W_J(x)),~x\in E\}\subseteq (-\infty,\infty]^{J+1}.
\end{eqnarray*}
The set ${\bf O}$ is convex in the cone $(-\infty,\infty]^{J+1}$ because the functions $W_j(\cdot)$ are affine and the space $E$ is convex. Since, by assumption, problem (\ref{eB12p}) is consistent and with a finite (optimal) value, ${\bf O}\cap \RR^{J+1}$ is a nonempty convex set in $\RR^{J+1}.$

Now we pass problem (\ref{eB12p}) to the following one in the space of performance vectors:
\begin{eqnarray}\label{eB02}
\mbox{Minimize } && W_0\\
\mbox{subject to } && \vec W=(W_0,W_1,\ldots,W_J)\in{\bf O}\nonumber\\
\mbox{and } && W_j\le d_j,~j=1,2,\ldots, J.\nonumber
\end{eqnarray}

The rest of the proof consists of verifying the statements formulated in each of the following steps.

\textit{Step 1.} There exists an optimal solution, say $\vec W^*=(W^*_0,W^*_1,\ldots,W^*_J)$, to problem (\ref{eB02}) such that
\begin{eqnarray*}\vec W^*=(W^*_0,W^*_1,\ldots,W^*_J)\in Par({\bf O}\cap\RR^{J+1}).
\end{eqnarray*}
(We mention that the assumption of totality of $\hat{{\cal C}}(E)$ is not in use here as well as in Steps 2 and 3 of this proof.)

To show the existence of such an optimal solution $\vec W^\ast$ to problem (\ref{eB02}), we consider the following problem
\begin{eqnarray}\label{PZ23aEqn01}
\mbox{Minimize over $x\in E:$} && \sum_{j=0}^J W_j(x)\\
\mbox{subject to } && W_j(x)\le d_j,~j=0,1,2,\ldots, J. \nonumber
\end{eqnarray}
Since $W_j(\cdot)$ is affine and lower semicontinuous, bounded from below, and with values in $(-\infty,\infty]$, so is the function $\sum_{j=0}^{J}W_j(\cdot)$.
Since $E$ is compact, the set of feasible solutions for problem (\ref{PZ23aEqn01}) is compact, as a closed subset of $E$. It follows that problem (\ref{PZ23aEqn01}) has an optimal solution $x^\ast\in E$.

Upon passing problem (\ref{PZ23aEqn01}) to the following problem in the space of performance vectors
\begin{eqnarray}\label{PZ23AEqn02}
\mbox{Minimize } && \sum_{j=0}^JW_j\\
\mbox{subject to } && \vec W=(W_0,W_1,\dots,W_J)\in{\bf O}\nonumber\\
\mbox{and } && W_j\le d_j,~j=0,1,2,\dots, J,\nonumber
\end{eqnarray}
we see that $\vec W^\ast=(W^*_0,W^*_1,\dots,W^*_J):=(W_0(x^\ast),W_1(x^\ast),\dots,W_J(x^\ast))$ is an optimal solution to problem (\ref{PZ23AEqn02}).

We argue that $\vec W^\ast$ is as required in Step 1 as follows.  Note that, compared to problem  (\ref{eB02}), in problem (\ref{PZ23aEqn01}), there is an additional constraint $W_0(x)\le d_0$, where $d_0$ is defined as the optimal value of problem (\ref{eB02}), which is finite by assumption, see the beginning of this proof. This implies that problem (\ref{PZ23aEqn01}) has a finite optimal value, so that $\vec{W}^\ast\in\RR^{J+1}$, and the set of feasible solutions for problem (\ref{PZ23aEqn01}) is a subset of the one for problem   (\ref{eB02}).  Therefore, $x^\ast$ is a feasible solution for problem (\ref{eB12p}), and accordingly, $\vec W^\ast$ is a feasible solution to problem  (\ref{eB02}).  Secondly, it necessarily holds that $W^\ast_0=d_0$, because otherwise, we would have $W_0(x^\ast)<d_0$, and this is against that $d_0$ is the optimal value of problem (\ref{eB12p}) and problem  (\ref{eB02}). Consequently, $\vec W^\ast$ is an optimal solution to problem  (\ref{eB02}), as required.

To see why $\vec W^\ast\in Par({\bf O}\cap \mathbb{R}^{J+1})$, note that $\vec W^\ast\in \mathbb{R}^{J+1}$, and if there is some $\vec{W}'=(W'_0,\dots,W'_J)\in {\bf O}\cap \mathbb{R}^{J+1}$, which strictly outperforms $\vec W^\ast$, i.e., $W'\le \vec W^\ast$ and $W'_j< W^\ast_j$ for some $j\in\{0,\dots,J\}$, then $\vec{W}'$ is feasible for problem  (\ref{PZ23AEqn02}) and $\sum_{j=0}^J W'_j<\sum_{j=0}^J W^\ast_j$. Consequently,
 $\vec W^\ast$ would not be optimal for problem  (\ref{PZ23AEqn02}), yielding a contradiction.  Step 1 is completed.

In Step 2, we consider $F$, a nonempty convex closed (and thus compact) subset of $E$.  Introduce the set
\begin{eqnarray}\label{TonightPZ23aEqn02}
\tilde{\bf O}:=\{\vec W(x):~x\in F\}\subseteq (-\infty,\infty]^{J+1},
\end{eqnarray}
where $\vec W(x)=(W_0(x),\dots,W_J(x))$.
The set $\tilde{\bf O}$ is a convex subset of $\bf O$, because $W_j(\cdot)$ are affine on $E$ and $F$ is a convex subset of $E$.

We assume that $\beta\in\RR$ and a vector $\vec b=(b_0,\dots,b_J)>0$ in $\RR^{J+1}$ are such that, for some $\vec u=(u_0,\dots,u_J)\in\tilde{\bf O}$, $\sum_{j=0}^J b_j u_j=\beta$ and
\begin{eqnarray}\label{May3PZ23aEqn055}
\sum_{j=0}^J b_j W_j \ge \beta~\forall~\vec W=(W_0,\dots,W_J)\in\tilde{\bf O}\cap\RR^{J+1}.
\end{eqnarray}

Recall that $0\times\infty:=0$ and $a+\infty:=\infty$ for all $a\in\RR\cup\{\infty\}$. All the inequalities for vectors are component-wise. Then from $\vec b>0$, we see that
\begin{eqnarray}\label{TonightPZ23aEqn01}
&&\hat{{\bf O}}:=\left\{\vec W=(W_0,\dots,W_J)\in\tilde{\bf O}:~\sum_{j=0}^J b_jW_j=\beta\right\}\nonumber\\
&=&
\left\{\vec W=(W_0,\dots,W_J)\in\tilde{\bf O}\cap \RR^{J+1}:~\sum_{j=0}^J b_jW_j\le\beta\right\}\subseteq \mathbb{R}^{J+1}.
\end{eqnarray}


\textit{Step 2.}  For $F$, $\vec{b}$ and $\beta$ as described in the above, the set
$\hat{{\bf O}}$
is bounded and closed, and thus is a compact subset of $\mathbb{R}^{J+1}$.

The set $\hat{\bf O}$ is bounded from below because each of the functions $W_j(\cdot)$ is bounded below by a constant say $\underline{W}_j$, and $\hat{\bf O}$ is a subset of the space of performance vectors ${\bf O}$.

The set $\hat{\bf O}$ is bounded from above because for each $j=0,1,\dots,J$ and $\vec W=(W_0,\dots,W_J)\in \hat{\bf O}$, $W_j= \frac{\beta-\sum_{i\ne j}b_i W_i}{b_j}\le \frac{\beta-\underline{W}_j\sum_{i\ne j}b_i}{b_j}<\infty.$ Again, here we used the fact that ${\vec b}>0$ componentwise.

The rest verifies the closedness of $\hat{\bf O}\subseteq \RR^{J+1}$ in $\RR^{J+1}.$

Let $\{\vec{W}^{(n)}=(W^{(n)}_0,\dots,W^{(n)}_J)\}_{n\ge 0}\subseteq \hat{\bf O}$ converge to some $\vec{W}=(W_0,\dots,W_J)\in \mathbb{R}^{J+1}$, i.e., $W_j^{(n)}\rightarrow W_j$ as $n\rightarrow \infty$ for each $j=0,\dots,J.$
Then
\begin{eqnarray}\label{3MayPZ23aEqn01}
\sum_{j=0}^J b_jW_j=\beta.
\end{eqnarray} It remains to show that $\vec{W}=(W_0(x),\dots,W_J(x))=:\vec{W}(x)$ for some $x\in F.$

For each $n\ge 0$, since, $\vec{W}^{(n)}\in \hat{\bf O}$, there is some $x_n\in F$ such that $\vec{W}^{(n)}=\vec{W}(x_n)$.

Since $F$ is a compact subset of $E$, by \cite[Theorem 2.31]{Aliprantis:2007}, the sequence $\{x_n\}_{n\ge 0}$ has a convergent subnet $\{y_\lambda\}_{\lambda\in\Lambda}$ with $y_\lambda\rightarrow x $ for some $x\in F$. We verify that $\vec{W}(x)=\vec{W}$ as follows.


Since for each $j=0,\dots,J,$ the function $W_j(\cdot)$ is lower semicontinuous,  by \cite[Lemma 2.42]{Aliprantis:2007},
\begin{eqnarray}\label{May3PZ23aEqn02}
\liminf_{\lambda}W_j(y_\lambda)\ge W_j(x)~ \forall~j=0,1,\dots,J.
\end{eqnarray}
 On the other hand, since for each $j\in\{0,1,\dots,J\},$ $\{W_j(y_\lambda)\}_{\lambda\in\Lambda}$ is a subnet of  $\{W_j(x_n)\}_{n\ge 0}=\{W_n^{(n)}\}_{n\ge 0}$ and
$W_j^{(n)}\rightarrow W_j$ as $n\rightarrow \infty$, by \cite[Lemma 2.17]{Aliprantis:2007}, $\{W_j(y_\lambda)\}_{\lambda\in\Lambda}\rightarrow W_j$, as its unique limit. It follows from this and (\ref{May3PZ23aEqn02}) that
\begin{eqnarray*}
W_j\ge W_j(x)
\end{eqnarray*}
for all $j=0,1\dots,J.$ In particular, $\vec{W}(x)\in\mathbb{R}^{J+1}.$

For the desired relation $\vec W= \vec W(x)$, it remains to show for each $j=0,1,\dots,J$ that the above inequality cannot hold strictly.
Suppose for contradiction that $W_j>W_j(x)$ for some $j\in\{0,1\dots,J\}$. Then we would have $\beta=\sum_{j=0}^J b_j W_j>\sum_{j=0}^J b_j W_j(x)\ge \beta$, where the first equality is by (\ref{3MayPZ23aEqn01}), the second inequality holds because  $\vec b>0$, and the last inequality holds by (\ref{May3PZ23aEqn055}) and the facts that $x\in F$, $\vec W(x)\in\mathbb{R}^{J+1}$, and $\tilde{\bf O}$ is the set of performance vectors of the points in $ F$. Hence, $\vec{W}=\vec{W}(x)$, as needed, and Step 2 is completed.

\textit{Step 3.} Consider the point $\vec{W}^\ast=(W^\ast_0,\dots,W^\ast_J)\in Par({\bf O}\cap \RR^{J+1})$ coming from Step 1. Then $\vec W^*=\sum_{k=1}^{J+1} \alpha_k \vec W_k$, where $\alpha_k\in[0,1]$ for each $k=1,\dots,J+1$, $\sum_{k=1}^{J+1} \alpha_k=1$, and $\vec W_k$ is an extreme point of ${\bf O}\cap \RR^{J+1}$, satisfying $\vec W_k\in Par({\bf O}\cap \RR^{J+1})$.

The justification of the claimed result in Step 3 is as follows.

Since $\vec{W}^\ast\in Par({\bf O}\cap \RR^{J+1})$ as asserted in Step 1, and ${\bf O}\cap \RR^{J+1}$ is a nonempty convex set in $\RR^{J+1}$ as was noted in the beginning of this proof,  Proposition \ref{lb31} is applicable, from which we deduce the following:
\begin{itemize}
\item[(a)] $G_{{\bf O}\cap \RR^{J+1}}(\vec W^\ast)\subseteq Par({\bf O}\cap \RR^{J+1})$.
\item[(b)] For some $1\le k\le J+1$, and $1\le i\le k$ there exist some $\vec{b}^i=(b^i_0,b^i_1,\dots,b^i_{J})\in \mathbb{R}^{J+1}$ and $\beta^i\in\mathbb{R}$, defining the hyperplanes \begin{eqnarray*}
H^i=\left\{\vec{p}=(p_0,\dots,p_J)\in\mathbb{R}^{J+1}:~\sum_{j=0}^{J} b^i_j p_j =\beta^i\right\},~i=1,2,\ldots, k,
\end{eqnarray*}
satisfying the following properties:
\begin{itemize}
\item[(i)] $\vec{b}^i\ge 0$, $\sum_{j=0}^J b^i_j> 0$, for $i=1,2,\ldots, k-1$ and $\vec{b}^k>0$. Here all the inequalities between vectors are componentwise.
\item[(ii)] $\sum_{j=0}^{J} b_j^1 p_j\ge \beta^1$ for all $\vec{p}=(p_0,\dots,p_J)\in {\bf O}\cap \RR^{J+1}$, $\sum_{j=0}^{J} b_j^1 W^\ast_j=\beta^1$; for $i=1,2,\dots,k-1$,   $\sum_{j=0}^{J} b_j^{i+1} p_j\ge \beta^{i+1}$ for all $\vec{p}=(p_0,\dots,p_J)\in {\bf O}\cap \RR^{J+1} \cap H^1\cap H^2\cap\dots\cap H^i $, and $\sum_{j=0}^{J} b_j^{i+1} W^\ast_j=\beta^{i+1}$.
\item[(iii)] $
G_{{\bf O}\cap \RR^{J+1}}(\vec W^\ast)= {\bf O}\cap \RR^{J+1} \cap H^1\cap H^2\cap\dots\cap H^k.
$
\end{itemize}
\end{itemize}
In the rest of this proof, we shall refer to the above consequences of  Proposition \ref{lb31} as the Feinberg-Shwartz lemma.

Let us verify that the set $G_{{\bf O}\cap \RR^{J+1}}(\vec W^\ast)$ is a convex and compact subset of $\mathbb{R}^{J+1}.$ First, we show its compactness by applying the statement established in Step 2.
To this end, we shall identify the appropriate $F,\vec{b}$ and $\beta$ satisfying the conditions described above Step 2.



Consider the half-spaces
\begin{eqnarray*}
\overline{H}^i:=\left\{\vec{p}=(p_0,\dots,p_{J})\in\mathbb{R}^{J+1}:~\sum_{j=0}^{J} b^i_j p_j \le \beta^i\right\},~i=1,2,\ldots, k.
\end{eqnarray*}
From Items (ii,iii) of (b) in the Feinberg-Schwartz lemma, we see that
\begin{eqnarray*}
&&G_{{\bf O}\cap \RR^{J+1}}(\vec W^\ast)= {\bf O}\cap \RR^{J+1} \cap H^1\cap H^2\cap\dots\cap H^k\\
&=&  \RR^{J+1} \cap {\bf O}\cap  \overline{H}^1\cap \overline{H}^2\cap\dots \cap \overline{H}^{k}
\end{eqnarray*}

Since the bounded from below functions $W_j(\cdot)$, $j=0,1,\dots,J$, are lower semicontinuous and $b_j^i\ge 0$ for all $j=0,1,\dots,J$, we see that $\sum_{j=0}^J b_j^i W_j(\cdot)$ is also lower semicontinuous. Here $0\cdot \infty:=0$ was in use.
It follows that the set
\begin{eqnarray*}
&&\{x\in E: \vec W(x) \in {\bf O}\cap \overline{H}^1\cap \overline{H}^2\cap\dots\cap \overline{H}^{k-1}\}\\
&=&\{x\in E: \vec W(x) \in  \overline{H}^1\cap \overline{H}^2\cap\dots\cap \overline{H}^{k-1}\}\\
&=&\bigcap_{i=1}^{k-1}\left\{x\in E: \sum_{j=0}^J b_j^i W_j(x)\le \beta^{i}\right\}.
\end{eqnarray*}
is closed and thus compact in $E$ because $E$ is compact.

If we take $F=\{x\in E: \vec W(x) \in {\bf O}\cap \overline{H}^1\cap \overline{H}^2\cap\dots\cap \overline{H}^{k-1}\}$, then the set in (\ref{TonightPZ23aEqn02}) takes the form $\tilde{\bf O}={\bf O}\cap \overline{H}^1\cap \overline{H}^2\cap\dots\cap \overline{H}^{k-1}$,  and
\begin{eqnarray*}
G_{{\bf O}\cap \RR^{J+1}}(\vec W^\ast)=\left\{\vec W=(W_0,\dots,W_J)\in\tilde{\bf O}\cap \RR^{J+1}:~\sum_{j=0}^J b_j^kW_j\le\beta^k\right\},
\end{eqnarray*}
which is in the same form as the set $\hat{{\bf O}}$ defined in (\ref{TonightPZ23aEqn01}) above Step 2. The conditions on $\vec u$, $b_j$ and $\beta$ above Step 2 are satisfied by $\vec{W}^\ast$, $b_j^k$ and $\beta^k$ by (b) in the Feinberg-Shwartz lemma. Hence, the statement in Step 2 is applicable, from which we see that  the set $G_{{\bf O}\cap \RR^{J+1}}(\vec W^\ast)$ is a compact subset in $\mathbb{R}^{J+1}$.  Its convexity is evident, and  is actually  by definition.

Now, since $G_{{\bf O}\cap \RR^{J+1}}(\vec W^\ast)\subseteq \mathbb{R}^{J+1}$ is an intersection with some hyperplane(s), its dimension is not higher than $J.$ One may apply the Caratheodory theorem together with the Krein-Milman theorem, see e.g., \cite[Theorems 5.32, 7.68]{Aliprantis:2007} or \cite[Corollary B.2.1]{APZY}, for that every point of  $G_{{\bf O}\cap \RR^{J+1}}(\vec W^\ast)$ can be written as the convex combination of at most $J+1$ extreme points in it. Since, according to (a) in the Feinberg-Shwartz lemma, $G_{{\bf O}\cap \RR^{J+1}}(\vec W^\ast)\subseteq Par({\bf O}\cap \RR^{J+1})$, we see that each of these extreme points are in $Par({\bf O}\cap \RR^{J+1})$. Since $G_{{\bf O}\cap \RR^{J+1}}(\vec W^\ast)$ is a face of ${\bf O}\cap \RR^{J+1}$, its extreme points are also extreme points of ${\bf O} \cap \RR^{J+1}$, see Lemma \ref{PZ23aRem01}. The statement in Step 3 is thus justified.

\textit{Step 4.} Consider an  extreme point $\vec W=(W_0,\dots,W_J)$ of ${\bf O}\cap\RR^{J+1}$, and suppose ${\vec W}\in Par({\bf O}\cap\RR^{J+1})$. Then there exists some extreme point $\hat{x}\in E$ such that $\vec W=\vec W(\hat x)=(W_0(\hat{x}),\dots,W_J(\hat{x}))$.

We justify the statement in Step 4 as follows.

Consider the set $E':=\{x\in E:~\vec W(x)=\vec W\}$.
Then $E'\ne\emptyset$ as $\vec{W}\in Par({\bf O}\cap\RR^{J+1})\subseteq {\bf O}$ by assumption, and
\begin{eqnarray}\label{TonightTomorrowPZ23aEqn01}
E'=\bigcap_{j=0}^J \{x\in E:~W_j(x)\le W_j\}.
\end{eqnarray}
Indeed, $E'$ is a subset of the set on the right-hand side. Let us justify that the opposite inclusion holds. Suppose that $x\in{E}$ belongs to the set on the right-hand side, i.e., it is such that $W_j(x)\le W_j$, $j=0,1,\ldots, J$. Then $\vec W(x)=(W_0(x,\dots,W_J(x)))\in{\bf O}\cap \RR^{J+1}$ and $\vec W(x)=\vec W$, because $\vec W\in Par({\bf O}\cap\RR^{J+1})$. This shows that the given point $x$ belongs to $E'$, as required.

Next, observe that $E'$ is a nonempty compact and convex set in $E$. Indeed, the nonemptiness of $E'$ was observed earlier. Since all the functions $W_j(\cdot)$, $j=0,1,\ldots, J$, are lower semicontinuous, it follows from (\ref{TonightTomorrowPZ23aEqn01}) that the set $E'$ is a closed and hence compact subset of the compact set $E$. It is obviously convex, as the functions $W_j(\cdot)$, $j=0,1,\ldots, J,$ are affine.

Applying Lemma \ref{lB6}(b) to the nonempty compact and convex set $E'$, we see that $E'$ contains some extreme point, say ${\hat x}\in E'$. Now we show that ${\hat x\in E'}$ is also an extreme point of  $E.$ By the definition of the set $E'$, this would lead to the statement claimed in Step 4.

Let $x_1,x_2\in E$ be such that $\hat{x}=\alpha x_1+(1-\alpha)x_2$ for some $\alpha\in (0,1)$.  Then
\begin{eqnarray}\label{GotoBedEqn01}
\vec{W}=\vec{W}(\hat{x})=\alpha \vec{W}(x_1)+(1-\alpha)\vec{W}(x_2),
\end{eqnarray}
where the first equality holds because $\hat{x}\in E'$, and the second equality holds as the functions $W_j(\cdot)$ are affine. Since $\vec{W}\in Par({\bf O}\cap \RR^{J+1})\subseteq \RR^{J+1}$, it follows that $\vec{W}(x_1),\vec{W}(x_2)\in\RR^{J+1}.$

If $ \vec{W}(x_1)= \vec{W}(x_2)$, then their common value is necessarily $\vec{W}$. Consequently, $x_1,x_2\in E'$. Since $\hat{x}$ is an extreme point of $E'$, we see that $x_1=x_2=\hat{x}$.

On the other hand,   it cannot happen that $ \vec{W}(x_1), \vec{W}(x_2)$ do not coincide, for otherwise, since  $\vec{W}(x_1),\vec{W}(x_2)\in {\bf O}\cap \mathbb{R}^{J+1}$,  (\ref{GotoBedEqn01}) indicates that $\vec W$ is not an extreme point of ${\bf O}\cap \RR^{J+1}$, which is a contradiction.

Therefore, $\hat{x}$ is also an extreme point of $E$, as requested. Step 4 is completed.

Now by Steps 3 and 4, we see that $\vec{W}^\ast$ coming from Step 1 satisfies $\vec{W}^\ast=\sum_{k=1}^{J+1}\alpha_k \vec{W}_k=\sum_{k=1}^{J+1}\alpha_k \vec{W}_k(x_k)$ for some $\alpha_k\in[0,1]$ satisfying $\sum_{k=1}^{J+1}\alpha_k=1$ and $x_k$ being an extreme point of $E$ for each $k=1,\dots,J+1.$ Since the functions $W_j(\cdot)$ are affine and $E$ is convex, from the previous equalities, we see $\vec{W}^\ast=\vec{W}(\sum_{k=1}^{J+1}\alpha_k x_k)$. Since $\vec{W}^\ast$ is an optimal solution for problem (\ref{eB02}),  the last equality shows that the point $\sum_{k=1}^{J+1}\alpha_k x_k\in E$ is an optimal solution for problem (\ref{eB12p}). It exhibits all the properties stated in Theorem \ref{tB1}.

Now, we have seen that the statement of this theorem holds, given the extra assumption that there is a feasible solution $\hat{x}\in E$ with  $W_0(\hat x)<\infty$.

Finally, we show that the theorem still holds, when the aforementioned assumption does not hold. Thus, until the end of this proof, we suppose that, for all feasible solutions $x\in E$, $W_0(x)=+\infty$.

If $J=0$, then any extreme point of $E$ is optimal. By Lemma \ref{lB6}(b), $E$ has at least one extreme point, and the statement of the theorem follows.

Consider the case of $J\ge 1$. In this case, we consider a modification  of problem (\ref{eB12p}):
\begin{eqnarray*}
&&\mbox{Minimize over $x\in E$:~} W_1(x)\\
&&\mbox{subject to~} W_j(x)\le d_j,~j=2,\dots,J.
\end{eqnarray*}
Since problem (\ref{eB12p}) is consistent by assumption,  an optimal solution to the above problem is necessarily feasible for problem (\ref{eB12p}), and its optimal value is finite and smaller or equal to $d_1$. Note also that this modified problem has $J-1$ constraints.  Now, we may apply what was proved earlier to this modified problem, and obtain an optimal solution for the modified problem in the form of a convex combination of at most $J$ extreme points of $E$. As was mentioned earlier, this solution is also optimal for the original problem (\ref{eB12p}).
This theorem is thus proved.
$\hfill\Box$

\section{Conclusion}\label{PZ23aSec06}
In conclusion, we considered an optimization problem in a convex set of a cone $E$ with an affine objective and $J$ affine constraints. The set $E$ is not required to be embedded in any vector space. Under suitable conditions, by applying the Feinberg-Shwartz lemma in finite dimensional convex analysis, we showed that there exists an optimal solution in the form of a mixture of at most $J+1$ extreme points of $E$. The result will be used in the study of Markov decision processes in a subsequent paper.


\end{document}